\documentclass[11pt,reqno]{amsart}

\usepackage[T1]{fontenc}
\usepackage[utf8]{inputenc}
\usepackage{lmodern}
\usepackage[english]{babel}
\usepackage{amsmath,amssymb,amsthm}
\usepackage{mathtools}
\usepackage{microtype}
\usepackage[a4paper,margin=30mm]{geometry}
\usepackage{enumitem}
\usepackage[hidelinks]{hyperref}

\newtheorem{theorem}{Theorem}[section]
\newtheorem{lemma}[theorem]{Lemma}
\newtheorem{proposition}[theorem]{Proposition}
\newtheorem{corollary}[theorem]{Corollary}
\newtheorem{maintheorem}{Theorem}

\theoremstyle{definition}
\newtheorem{definition}[theorem]{Definition}
\theoremstyle{remark}

\numberwithin{equation}{section}

\DeclareMathOperator{\can}{can}
\newcommand{\Can}{\mathrm{Can}}
\newcommand{\Rel}{\mathrm{Rel}}
\newcommand{\Z}{\mathbb Z}
\newcommand{\ncl}[1]{\langle\!\langle #1\rangle\!\rangle}

\title[Nonamenability for free Burnside groups]{A new proof of nonamenability for free Burnside groups of odd exponent}

\author{Anna Cascioli}
\email[Anna Cascioli]{acasciol@uni-muenster.de}

\author{Sam Shepherd}
\email[Sam Shepherd]{sshepher@uni-muenster.de}

\author{Katrin Tent}
\email[Katrin Tent]{tent@uni-muenster.de}

\address{University of M{\"u}nster, Einsteinstrasse 62, M{\"u}nster 48149, Germany}

\date{}
\hypersetup{
  pdftitle={A new proof of nonamenability for free Burnside groups of odd exponent},
  pdfauthor={Anna Cascioli, Sam Shepherd, Katrin Tent}
}

\begin{document}

\begin{abstract}
We prove that the free Burnside group $B(m,n)$ is nonamenable for every $m\ge 2$ and every odd exponent $n\ge 557$. This improves Adian's bound $n\ge 665$ and matches the infiniteness bound of Atkarskaya, Rips and Tent.
Our proof combines their canonical form construction with a counting argument to obtain an explicit isoperimetric inequality.
\end{abstract}

\maketitle

\section{Introduction}\label{sec:introduction}

In 1902, Burnside asked whether a finitely generated group must be finite if all its elements have finite order \cite{Burnside1902}.
The bounded Burnside problem asks the same question under the stronger assumption that the group satisfies the law $x^n=1$ for some fixed positive integer $n$. Every $m$-generated group satisfying this law is a quotient of the free Burnside group
\[
  B(m,n)=F_m/\ncl{w^n\mid w\in F_m},
\]
where $F_m$ is the free group of rank $m$.

The bounded problem was answered negatively in 1968, when Novikov and Adian proved that $B(m,n)$ is infinite whenever $m\ge 2$ and $n\ge 4381$ is odd
\cite{NovikovAdian1968I,NovikovAdian1968II,NovikovAdian1968III}.
Adian subsequently lowered the bound to $n\ge 665$ \cite{Adian1979}. In 2023, Atkarskaya, Rips and Tent proved infiniteness for every odd $n\ge 557$ \cite{ART}, giving the best currently known bound for odd exponents.

Beyond infiniteness, a central question concerns amenability. Free Burnside groups are torsion groups and therefore contain no nonabelian free subgroups, so their nonamenability must be established by other means. Adian achieved this in 1982 for $m\ge 2$ and odd $n\ge 665$ \cite{Adian1983}, using estimates on the growth of freely reduced words representing the identity together with Grigorchuk's cogrowth criterion \cite{Grigorchuk1980}. Using the canonical form construction of Atkarskaya, Rips and Tent,
we prove the following quantitative form of nonamenability.

\begin{maintheorem}\label{thm:main}
Let $m\ge 2$, let $n\ge 557$ be odd, and let $a,b$ be two distinct standard generators of $B(m,n)$. Then every finite subset $F\subseteq B(m,n)$ satisfies
\begin{equation}\label{eq:main}
  |aF\setminus F|+|bF\setminus F|
  \ge \frac{2}{5}\left(1-\frac4n\right)|F|.
\end{equation}
In particular, $B(m,n)$ is nonamenable.
\end{maintheorem}

This lowers Adian's nonamenability bound from $665$ \cite{Adian1983} to $557$, matching the best currently known infiniteness bound for odd exponents \cite{ART}. The boundary estimate is independent of the rank and involves only two standard generators.

Our proof approaches nonamenability through the boundaries of finite subsets. The canonical form construction in \cite{ART} uses iterated small cancellation to produce canonical representatives by induction on the nesting depth of relators. More precisely, for the quotient map $\pi:F_m\to B(m,n)$, it yields a  section
\[
  \can:B(m,n)\longrightarrow F_m,
  \qquad \pi\circ\can=\operatorname{id}_{B(m,n)}.
\]
Thus each element of $B(m,n)$ is assigned a distinguished reduced word, its \emph{canonical form}, and equality of group elements can be tested by literal equality of their canonical representatives. We recall the properties of these representatives needed here in Section~\ref{sec:preliminaries}. Our proof combines these properties with a counting argument that imposes no further restriction on the
exponent.

\subsection{Proof strategy}\label{subsec:proof-strategy}

Put $A=\langle a\rangle$. We first prove the boundary estimate for finite sets $F$ satisfying $aF=F$, which are unions of $A$-cosets, each containing $n$ elements. The $b$-boundary counts the points that leave $F$ under multiplication by $b$, so it suffices to bound the number of internal connections $(x,bx)$ with both endpoints in $F$.

To control these connections, we consider chains that move within an $A$-coset by a balanced power of $a$ and then cross a $b$-connection in either direction. We encode each chain by its successive exponents and directions, and require the encoding to avoid immediate free cancellation and consecutive repetitions of a nonempty block. For these admissible chains, the multiplier of every nonempty subchain is a canonical word containing $b^{\pm1}$ and hence represents an element
outside $A$. The chain therefore cannot revisit an $A$-coset. Yet if every coset had at least five incident connections, a counting argument would produce an admissible chain visiting more cosets than $F$ contains. It follows that some coset has at most four connections; removing such cosets successively yields the required $b$-boundary estimate.

Finally, we pass to an arbitrary finite set $F$ by filling every $A$-coset that meets it. There are at most $|F|/n$ complete cosets, and each incomplete coset contributes to the $a$-boundary. Combining this coset count with the invariant estimate gives a weighted boundary inequality. Interchanging $a$ and $b$ and adding the two resulting inequalities yields the estimate in Theorem~\ref{thm:main}.

Section~\ref{sec:preliminaries} recalls the F\o lner criterion and the canonical form properties needed below. Section~\ref{sec:chains} introduces internal connections and admissible chains, and proves that such chains cannot return to a coset. The counting argument is given in Section~\ref{sec:counting}. Section~\ref{sec:boundary} deduces the boundary estimates and proves Theorem~\ref{thm:main}.

\section{Preliminaries}\label{sec:preliminaries}

\subsection{Amenability and the F\o lner criterion}\label{subsec:folner}

A discrete group $G$ is \emph{amenable} if it admits a left-invariant mean on $\ell^\infty(G)$, that is, a positive normalized linear functional invariant under left translations. We use the following characterization in terms of finite subsets.

\begin{theorem}[F\o lner criterion {\cite{Folner1955}}]\label{thm:folner}
A discrete group $G$ is amenable if and only if, for every finite subset $S\subseteq G$ and every $\varepsilon>0$, there is a nonempty finite set $F\subseteq G$ such that
\[
  \sum_{s\in S}|sF\setminus F|<\varepsilon|F|.
\]
\end{theorem}

For a finite set $F$, one has $|sF\mathbin{\triangle}F|=2|sF\setminus F|$, so this formulation is equivalent to the usual one with symmetric differences. In particular,
to prove nonamenability it suffices to find a finite subset $S$ and a constant $c>0$ such that
\[
  \sum_{s\in S}|sF\setminus F|\ge c|F|
\]
for every finite subset $F$ of the group. We will establish such an estimate with $S=\{a,b\}$, where $a,b$ are two distinct standard generators of $B(m,n)$.

\subsection{Notation and canonical forms}

Throughout the paper, $m\ge 2$ is an integer and $n\ge 557$ is odd. Fix distinct free generators $a,b$ of $F_m$, and use the same letters for their images in $B(m,n)$. For a written word $W$, let $[W]$ denote the element it represents in $B(m,n)$, and let $|W|$ denote its word length. Equality of words means literal equality, while equality in $B(m,n)$ is expressed using brackets when a distinction is needed. Put
\[
  A=\langle a\rangle,
  \qquad
  D=\left\{-\frac{n-1}{2},\ldots,\frac{n-1}{2}\right\}.
\]
The natural map $B(m,n)\to(\Z/n\Z)^m$ shows that $a$ has order $n$ and $b\notin A$. In particular, $D$ gives a complete set of exponents for the elements of $A$. We use right cosets $Ax$, which are the orbits of left multiplication by $A$.

We recall the features of the canonical form construction in \cite{ART} needed below. The construction organizes the defining relations into pairwise disjoint families $\Rel_r$, $r\ge 1$, according to the nesting of power subwords, with nesting constant $\tau=15$. Every relator is $c^n$ for a nonempty cyclically reduced word $c$ that is not a proper power. The families are closed under inverses and cyclic shifts, and
\[
  \Rel_1=\{x^n\mid x\text{ is a free generator or its inverse}\}.
\]
Writing $N_0=\{1\}$ and $N_r=\ncl{\Rel_1\cup\cdots\cup\Rel_r}$, one has $B(m,n)=F_m/\bigcup_{r\ge 0}N_r$.

At rank $r$, the map $\can_r$ selects a unique representative of each coset modulo $N_r$. Its image is denoted by $\Can_r$, and $\Can_0$ is the set of freely reduced words. These sets satisfy
\[
  \Can_0\supseteq\Can_1\supseteq\Can_2\supseteq\cdots,
\]
and $\can_r(W)=W$ whenever $W\in\Can_r$.
For every word $W$, the sequence $\can_r(W)$ eventually stabilizes at the representative $\can([W])$ chosen by the section introduced above. We abbreviate this value to $\can(W):=\can([W])$. The fundamental uniqueness property is
\begin{equation}\label{eq:canonical-uniqueness}
  [U]=[V]\quad\Longleftrightarrow\quad \can(U)=\can(V).
\end{equation}
For these properties, see \cite[Section~3, Proposition~8.2, Corollary~8.7 and Theorem~8.1]{ART}.

A \emph{fractional power of rank $r$} is a nonempty subword $u$ of $c^k$, where $k\in\mathbb Z$ and $c^n\in\Rel_r$. Its measure with respect to the specified period $c$ is $|u|/|c|$.
We will use the following two criteria for a word to remain canonical.

\emph{Rank one.}
If $W$ is reduced and every maximal block in a single generator has length at most $(n-1)/2$, then $\can_1(W)=W$ and $W\in\Can_1$ \cite[Section~5.2]{ART}.

\emph{Higher ranks.}
If $r\ge 2$ and $W\in\Can_{r-1}$ contains no fractional power of rank $r$ of measure greater than $n/2-5\tau-2$, then $W\in\Can_r$. This is the property stated as IH~4 in \cite[Section~3]{ART} and established at every rank by the construction. For the exponents under consideration,
\begin{equation}\label{eq:threshold}
  \frac n2-5\tau-2=\frac n2-77>3.
\end{equation}

These criteria will allow us to show that certain reduced words belong to $\Can_r$ for every rank $r$, and hence are already canonical.

\subsection{A criterion allowing long powers of one generator}

An \emph{$a$-block} is a nonempty maximal consecutive sequence of $a$'s or of $a^{-1}$'s. A \emph{cube} is a literal subword $V^3$, where $V$ is a nonempty word. The following consequence of the preceding criteria is the sole canonical form input to the combinatorial proof.

\begin{lemma}\label{lem:canonical-criterion}
Let $W$ be a reduced word in $a^{\pm1},b^{\pm1}$. Suppose that every $a$-block has length at most $(n-1)/2$ and that $W$ contains no cube $V^3$ whose period $V$ contains a letter $b^{\pm1}$. Then $\can(W)=W$. If $W$ contains a letter $b^{\pm1}$, then $[W]\notin A$.
\end{lemma}

\begin{proof}
Since $W$ is reduced, $W\in\Can_0$. Every $a$-block has length at most $(n-1)/2$, while every $b$-block has length at most two: otherwise, $W$ would contain $b^3$ or $b^{-3}$, contrary to the hypothesis. Thus all maximal single-generator blocks have length
at most $(n-1)/2$. The rank-one criterion gives $\can_1(W)=W$ and $W\in\Can_1$.

We now prove that $W$ remains canonical at every higher rank.
Suppose that $r\ge 2$ and, inductively, $W\in\Can_{r-1}$. By IH~4, recalled above, it suffices to exclude fractional powers of rank $r$ and measure greater than $n/2-5\tau-2$.
Suppose that $u$ is such an occurrence in $W$, with period $c$ satisfying $c^n\in\Rel_r$. By \eqref{eq:threshold},
\[
  \frac{|u|}{|c|}>\frac n2-5\tau-2>3.
\]
The first $3|c|$ letters of $u$ therefore form a literal cube $V^3$, where $V$ is a cyclic shift of $c$ or $c^{-1}$.

We claim that $V$ contains a letter $b^{\pm1}$. Since $V$ occurs in $W$, it uses only $a^{\pm1}$ and $b^{\pm1}$. If it contained no $b^{\pm1}$, reducedness would give $V=a^s$ for some nonzero integer $s$. Since $c$ is not a proper power, neither is its
cyclic shift $V$ or the corresponding cyclic shift of $c^{-1}$. Hence $|s|=1$. It follows that $|c|=|V|=1$, so $c^n\in\Rel_1$. This contradicts $c^n\in\Rel_r$, since $r\ge 2$ and the relator families are pairwise disjoint.

Thus $V$ contains $b^{\pm1}$, and the occurrence of $V^3$ in $W$ contradicts the hypothesis. There is therefore no such fractional power $u$, and IH~4 gives $W\in\Can_r$. By induction, $W$ belongs
to $\Can_r$ for every rank, so its stabilized canonical form is $\can(W)=W$.

Finally, every element of $A$ is represented by $a^d$ for some $d\in D$. Each such word satisfies the hypotheses of the lemma, so the first part of the proof gives $\can(a^d)=a^d$.
If $[W]\in A$, then $[W]=[a^d]$ for some $d\in D$, and uniqueness of canonical forms yields
\[
  W=\can(W)=\can(a^d)=a^d.
\]
This is impossible when $W$ contains a letter $b^{\pm1}$. Hence, $[W]\notin A$.
\end{proof}

Unlike an unrestricted cube-free condition, this criterion permits arbitrary balanced powers of $a$. This is essential below: moving within an $A$-coset may require any exponent in $D$.

\section{Admissible chains through cyclic cosets}\label{sec:chains}

\subsection{Internal connections}

Let $F\subseteq B(m,n)$ be a finite $a$-invariant set. Thus
\[
  F=Aw_1\sqcup\cdots\sqcup Aw_k,\qquad |F|=nk.
\]
The representatives $w_i$ are arbitrary elements of $B(m,n)$: no canonical representatives of the cosets are required. An \emph{internal $b$-connection} is an ordered pair $(x,bx)$ with both endpoints in $F$. Put
\[
  e(F)=\#\{x\in F\mid bx\in F\}=|bF\cap F|.
\]
Different pairs are counted separately even when they join the same two cosets. Since left multiplication by $b$ is a bijection,
\begin{equation}\label{eq:connections-boundary}
  |bF\setminus F|=nk-e(F).
\end{equation}
Thus an upper bound on internal connections gives a lower bound on the $b$-boundary.

For an $A$-coset $C\subseteq F$, define its number of incident connections by
\begin{equation}\label{eq:degree}
  d_F(C)=\#\{x\in C\mid bx\in F\} +\#\{x\in C\mid b^{-1}x\in F\}.
\end{equation}
The first term counts connections starting in $C$, and the second counts those ending in $C$. No connection has both endpoints in the same coset: if $x,bx\in C$, then $bx=a^h x$ for some integer $h$, so $b=a^h$, contrary to $b\notin A$. Consequently,
\begin{equation}\label{eq:deletion}
  e(F)=e(F\setminus C)+d_F(C).
\end{equation}
We will prove that every nonempty such $F$ has a coset $C$ with $d_F(C)\le 4$. Repeatedly deleting these cosets will then give $e(F)\le 4k$.

\subsection{Chains and their encodings}

Starting from $x\in F$, a step has the form
\[
  y=b^{\varepsilon}a^d x\in F,
  \qquad \varepsilon\in\{1,-1\},\quad d\in D.
\]
We first move within the current $A$-coset and then traverse a $b$-connection, forwards or backwards. The intermediate point $a^d x$ belongs to $F$ because $F$ is $a$-invariant. For every $x\in C$, the points $a^d x$, $d\in D$, enumerate $C$ exactly once, so
\begin{equation}\label{eq:available-pairs}
  \#\{(\varepsilon,d)\in\{1,-1\}\times D
        \mid b^{\varepsilon}a^d x\in F\}=d_F(C).
\end{equation}
Thus all connections incident to the current coset are available from any chosen point of that coset.

A \emph{chain of length $t$} consists of a starting point $x_0\in F$ and a list of pairs $(\varepsilon_1,d_1),\ldots,(\varepsilon_t,d_t)$ such that
\[
  x_j=b^{\varepsilon_j}a^{d_j}x_{j-1}\in F,
  \qquad 1\le j\le t.
\]
The exponent $d_j$ is measured from $x_{j-1}$, not from a fixed coset representative. The starting point and the pair list determine all points of the chain. We count distinct pair lists as distinct chains, even if they have the same endpoint.

For $0\le i<j\le t$, the written multiplier of the corresponding subchain is
\[
  U_{i,j}=b^{\varepsilon_j}a^{d_j}\cdots
             b^{\varepsilon_{i+1}}a^{d_{i+1}},
  \qquad x_j=[U_{i,j}]x_i.
\]
In particular,
\begin{equation}\label{eq:return}
  Ax_i=Ax_j\quad\Longrightarrow\quad [U_{i,j}]\in A.
\end{equation}
To prevent returns, we impose restrictions that allow Lemma~\ref{lem:canonical-criterion} to be applied to every nonempty subchain multiplier.

\begin{definition}\label{def:admissible}
A chain is \emph{admissible} if its pair list satisfies:
\begin{enumerate}[label=(\Alph*),ref=\Alph*,leftmargin=*]
  \item\label{cond:reduced}
  For $j\ge 2$, one never has both $d_j=0$ and
  $\varepsilon_j=-\varepsilon_{j-1}$.
  \item\label{cond:squarefree}
  The list is \emph{square-free}, i.e. it contains no contiguous sublist $QQ$, where $Q$ is a nonempty list of pairs.
\end{enumerate}
\end{definition}

Both conditions are inherited by subchains. Empty and one-step chains are admissible. Condition~(\ref{cond:reduced}) makes the multiplier freely reduced, while condition~(\ref{cond:squarefree})
prevents cubes involving $b$, as follows.

\begin{lemma}\label{lem:multiplier}
The multiplier word
\[
  U=b^{\varepsilon_t}a^{d_t}\cdots b^{\varepsilon_1}a^{d_1}
\]
of a nonempty admissible chain is reduced and canonical, and $[U]\notin A$.
\end{lemma}

\begin{proof}
Condition~(\ref{cond:reduced}) excludes every possible cancellation between adjacent factors, so $U$ is reduced. Every maximal $a$-block is one of the nonempty blocks $a^{d_j}$ and has length at most $(n-1)/2$. Reading the $b$-letters from left to right, together with the $a$-blocks following them, recovers the reversed pair list
\begin{equation}\label{eq:reversed-list}
  (\varepsilon_t,d_t),\ldots,(\varepsilon_1,d_1),
\end{equation}
which is also square-free.

Suppose that $U$ contains a cube $V^3$, where $V$ contains $\ell\ge 1$ letters $b^{\pm1}$. Number the $3\ell$ occurrences of these letters in the cube from left to right. For each of the first $2\ell$ occurrences, record its sign and the exponent of the $a$-block between it and the next $b$-letter, taking exponent zero for an empty gap. These are $2\ell$ consecutive complete pairs in \eqref{eq:reversed-list}: the final gap lies inside the cube because the $(2\ell+1)$-st
$b$-letter exists. Translation by $|V|$ letters carries the $j$-th $b$-letter to the $(j+\ell)$-th with the same sign and the same following gap, for $1\le j\le\ell$. Hence the first $\ell$ recorded
pairs repeat immediately, contradicting square-freeness.

Lemma~\ref{lem:canonical-criterion} now applies to $U$. Since $U$ contains a $b$-letter, it gives $\can(U)=U$ and $[U]\notin A$.
\end{proof}

\begin{corollary}\label{cor:no-return}
An admissible chain $x_0,\ldots,x_t$ visits $t+1$ distinct $A$-cosets. In a union of $k$ cosets, its length is therefore at most $k-1$.
\end{corollary}

\begin{proof}
If $Ax_i=Ax_j$ for some $i<j$, then \eqref{eq:return} gives $[U_{i,j}]\in A$. The subchain from $x_i$ to $x_j$ is nonempty and admissible, contradicting Lemma~\ref{lem:multiplier}.
\end{proof}

\section{Counting admissible chains}\label{sec:counting}

We now show that sufficiently many connections force a chain longer than Corollary~\ref{cor:no-return} permits. The counting step uses the entire pair list: when an extension creates a repeated terminal block, deleting its second copy gives a shorter admissible chain from which that extension can be reconstructed.

\begin{lemma}\label{lem:long-chain}
Suppose that $F$ is a union of $k\ge 1$ cosets of $A$ and that $d_F(C)\ge 5$ for every coset $C\subseteq F$. From every fixed point $x_0\in F$, there is an admissible chain of length $k$ inside $F$.
\end{lemma}

\begin{proof}
Fix a point $x_0\in F$. For $0\le t\le k$, let $\alpha_t$ be the number of admissible chains of length $t$ starting at $x_0$, with different pair lists counted separately. Each $\alpha_t$ is finite, since there are only $2n$ possible pairs at each step. The empty chain gives $\alpha_0=1$. We prove inductively that
\begin{equation}\label{eq:growth}
  \alpha_{t+1}\ge 2\alpha_t,\qquad 0\le t<k.
\end{equation}
This estimate will give the desired existence statement
$\alpha_k>0$.

By \eqref{eq:available-pairs} and the assumption $d_F(C)\ge 5$, there are at least five possible first steps. Every one-step chain is admissible, so $\alpha_1\ge 5\ge 2\alpha_0$.

Now let $1\le t<k$. At the endpoint of an admissible chain of length $t$, there are at least five available next pairs. At most one, namely $(-\varepsilon_t,0)$, is excluded by condition~(\ref{cond:reduced}). Consider all the remaining extensions and call them \emph{candidates}. There are at least four candidates per parent, and removing the last pair recovers the parent itself. Hence different parent lists give distinct candidates, and there are at least $4\alpha_t$ candidates in total.

The \emph{good} candidates are precisely the chains counted by $\alpha_{t+1}$. Every \emph{bad} candidate still satisfies condition~(\ref{cond:reduced}), but its pair list is no longer square-free. Since the parent was square-free, every square appearing in a bad candidate must end at the newly added last pair. Such a candidate therefore has the form $RQQ$, where $R$ is a possibly empty list and $Q$ has length
\[
  1\le q\le\left\lfloor\frac{t+1}{2}\right\rfloor.
\]

If we delete the second copy of $Q$, the retained list $u=RQ$ has length $t+1-q$ and is admissible. Indeed, it is a prefix of the admissible parent, since $q\ge 1$. For fixed $q$, this deletion is injective: $u$ has length at least $q$, its last $q$ pairs are $Q$, and appending those pairs reconstructs $RQQ$ uniquely. Since the starting point is fixed, the reconstructed pair list also determines all the group elements of the candidate. We only need uniqueness: an arbitrary shorter admissible chain need not reconstruct a valid candidate.

Therefore, at most $\alpha_{t+1-q}$ bad candidates have a terminal square with block length $q$. Summing over $q$ may count a candidate more than once, but this is harmless for the upper bound. We obtain
\begin{equation}\label{eq:recurrence}
  \alpha_{t+1}\ge
  4\alpha_t-
  \sum_{q=1}^{\lfloor(t+1)/2\rfloor}\alpha_{t+1-q}.
\end{equation}
Every index on the right is at most $t$, so these are counts from previously considered lengths.

Now assume that $\alpha_j\ge 2\alpha_{j-1}$ has been established for $1\le j\le t$. Then
\[
  \alpha_{t+1-q}\le 2^{-(q-1)}\alpha_t,
  \qquad
  1\le q\le\left\lfloor\frac{t+1}{2}\right\rfloor.
\]
For $q=1$ this is an equality; for $q=2$ it is $\alpha_{t-1}\le\alpha_t/2$, and so on. Hence the total number of bad candidates is at most
\[
  \alpha_t\left(1+\frac12+\frac14+\cdots\right)=2\alpha_t.
\]
Equation~\eqref{eq:recurrence} therefore gives
\[
  \alpha_{t+1}\ge 4\alpha_t-2\alpha_t=2\alpha_t,
\]
completing the induction. In particular, $\alpha_k\ge 2^k>0$, so an admissible chain of length $k$ exists.
\end{proof}

\section{Boundary estimates and nonamenability}\label{sec:boundary}

\subsection{The invariant case}

\begin{proposition}\label{prop:invariant}
Every nonempty finite union $F$ of $A$-cosets contains a coset $C$ with $d_F(C)\le 4$. If $F$ consists of $k\ge 0$ cosets, then
\begin{equation}\label{eq:invariant-estimate}
  |bF\cap F|=e(F)\le 4k,
  \qquad
  |bF\setminus F|\ge\left(1-\frac4n\right)|F|.
\end{equation}
\end{proposition}

\begin{proof}
If every coset had at least five incident connections, Lemma~\ref{lem:long-chain} would give an admissible chain of length $k$, contradicting Corollary~\ref{cor:no-return}. Hence some coset has at most four incident connections.

The bound $e(F)\le 4k$ follows by induction on $k$, with the case $k=0$ immediate. For $k\ge 1$, choose such a coset $C$. By \eqref{eq:deletion} and the induction hypothesis,
\[
  e(F)=e(F\setminus C)+d_F(C)\le 4(k-1)+4=4k.
\]
The boundary estimate follows from \eqref{eq:connections-boundary} and $|F|=nk$.
\end{proof}

\subsection{Arbitrary finite sets}

We conclude by comparing a finite set with the union of all $A$-cosets that meet it. The essential observation is that each incomplete coset contributes to the $a$-boundary.

\begin{proof}[Proof of Theorem~\ref{thm:main}]
The assertion is immediate for $F=\varnothing$, so assume that $F$ is nonempty. Its $A$-saturation
\[
  \widehat F=AF=\bigcup_{j=0}^{n-1}a^jF
\]
is finite and $a$-invariant. Let $k$ be its number of $A$-cosets. Since $F\subseteq\widehat F$, Proposition~\ref{prop:invariant} gives
\begin{equation}\label{eq:saturation-connections}
  |F|-|bF\setminus F|=|bF\cap F|
  \le |b\widehat F\cap\widehat F|\le 4k.
\end{equation}

Call a coset of $\widehat F$ \emph{complete} if it is contained in $F$, and \emph{incomplete} otherwise. There are at most $|F|/n$ complete cosets. If $C$ is incomplete, then $F\cap C$ is a nonempty proper subset of the $n$-cycle given by left multiplication by $a$. It follows that some $x\in F\cap C$ satisfies $ax\notin F$. Thus $C$ contains a point of $aF\setminus F$. Distinct incomplete cosets give distinct boundary points, so
\begin{equation}\label{eq:coset-count}
  k\le \frac{|F|}{n}+|aF\setminus F|.
\end{equation}
Combining \eqref{eq:saturation-connections} and \eqref{eq:coset-count}, we obtain
\begin{equation}\label{eq:weighted-boundary}
  4|aF\setminus F|+|bF\setminus F|
  \ge\left(1-\frac4n\right)|F|.
\end{equation}
The preceding argument applies to any ordered pair of distinct
standard generators. Interchanging $a$ and $b$ therefore gives
\[
  |aF\setminus F|+4|bF\setminus F|
  \ge \left(1-\frac4n\right)|F|.
\]
Adding this inequality to \eqref{eq:weighted-boundary}, we obtain
\[
  5\bigl(|aF\setminus F|+|bF\setminus F|\bigr)
  \ge 2\left(1-\frac4n\right)|F|.
\]
Dividing by $5$ proves \eqref{eq:main}.

Since $\frac{2}{5}(1-\frac4n)>0$, the F\o lner criterion (Theorem~\ref{thm:folner}) now gives nonamenability.
\end{proof}

\section*{Acknowledgments}
The first and second authors are funded by the Deutsche Forschungsgemeinschaft (DFG, German Research Foundation) under Germany's Excellence Strategy EXC 2044/2 –390685587, Mathematics Münster: Dynamics–Geometry–Structure.

\section*{AI Declaration}

All mathematical ideas, results, and proofs presented in this paper were developed without the assistance of large language models (LLMs).

\bibliographystyle{amsplain}
\bibliography{burnside_references}

\end{document}